\documentclass[reqno,12pt]{amsart}

\usepackage{amsmath}
\usepackage{amssymb}
\usepackage{amsthm}
\usepackage[english]{babel}
\usepackage{graphics, xcolor}
\newtheorem{theorem}{Theorem}

\newtheorem{lemma}{Lemma}

\newtheorem{corollary}{Corollary}
\newcommand{\beqa}{\begin{eqnarray}}
\newcommand{\beqan}{\begin{eqnarray*}}
\newcommand{\eeqa}{\end{eqnarray}}
\newcommand{\eeqan}{\end{eqnarray*}}
\def\beq#1\eeq{\begin{equation}#1\end{equation}}

 \def\na{\,\, {\raise.4pt\hbox{$\shortmid$}}{\hskip-2.0pt\to}\, \, }

\def\={\overset{ \text{\rm def} }=}

\newcommand{\bl}{}

\newtheorem{remark}{Remark}

\begin{document}

\title[Concentration Functions]
{New applications of Arak's inequalities\\ to  the
Littlewood--Offord problem}

\author[F.~G\"otze]{Friedrich G\"otze}
\author[A.Yu. Zaitsev]{Andrei Yu. Zaitsev}

\email{goetze@math.uni-bielefeld.de}
\address{Fakult\"at f\"ur Mathematik,\newline\indent
Universit\"at Bielefeld, Postfach 100131,\newline\indent Bielefeld, D-33501,
Germany\bigskip}
\email{zaitsev@pdmi.ras.ru}
\address{St.~Petersburg Department of Steklov Mathematical Institute
\newline\indent
Fontanka 27,
\newline\indent
St.~Petersburg, 191023, Russia
\newline\indent
and
\newline\indent
St.Petersburg State University,
\newline\indent
Universitetskaya nab. 7/9,
\newline\indent
St. Petersburg, 199034, Russia}

\begin{abstract}
Let $X_1,\ldots,X_n$ be independent identically distributed random variables.
In this paper we study the behavior of concentration functions of
 weighted sums $\sum_{k=1}^{n}X_ka_k $ with respect to the
arithmetic structure of coefficients~$a_k$ in the context of the
Littlewood--Offord problem. In  recent
 papers of Eliseeva, G\"otze and  Zaitsev, we
discussed the relations between the inverse principles stated by
Nguyen, Tao and Vu and similar principles formulated  by Arak in
his papers from the 1980's. In this paper, we will
derive some more general and more precise consequences
 of Arak's inequalities
providing new results in  the context of the
Littlewood--Offord problem.
\end{abstract}

\footnotetext[1]{The authors were supported by the SFB 701 in
Bielefeld and by SPbGU-DFG grant
6.65.37.2017.}\footnotetext[2]{The  second author was supported by
grant RFBR 16-01-00367 and by the Program of the Presidium of the Russian
Academy of Sciences No. 01 'Fundamental Mathematics and its Applications' under
grant PRAS-18-01.}

\keywords {concentration functions, inequalities, the
Littlewood--Offord problem, sums of independent random variables}

\subjclass {Primary 60G50; secondary 11P70, 60E07, 60E10, 60E15}

\maketitle

\section{Introduction}

The
 concentration function of an $\mathbf{R}^d$-{valued}
vector $Y$ with distribution $F=\mathcal L(Y)$ {is defined by}
\begin{equation}
Q(F,\tau)=\sup_{x\in\mathbf{R}^d}\mathbf{P}(Y\in x+ \tau B), \quad
\tau\geq0, \nonumber
\end{equation}
where $B=\{x\in\mathbf{R}^d:\|x\|\leq 1/2\}$ 
denotes the centered
Euclidean ball of radius 1/2.

Let $X,X_1,\ldots,X_n$ be independent identically distributed
(i.i.d.{}) random variables. Let $a=(a_1,\ldots,a_n)\ne 0$, where
$a_k=(a_{k1},\ldots,a_{kd})\in \mathbf{R}^d$, $k=1,\ldots, n$.
Starting with seminal papers of Littlewood and Offord \cite{LO}
and Erd\"os~\cite{Erd}, the behavior of the concentration
functions of the weighted sums $S_a=\sum\limits_{k=1}^{n} X_k a_k$
has been intensively  studied. Denote by $F_a=\mathcal L(S_a)$ the
distribution of~$S_a$. We refer to~\cite{EGZ2} for a discussion of
the history of the problem.

Several years ago, Tao and Vu \cite{Tao and Vu} and Nguyen and Vu
\cite{Nguyen and Vu} proposed the so-called 'inverse principles' in
the Littlewood--Offord problem (see  Section~\ref{2}). In the
 papers of G\"otze, Eliseeva and  Zaitsev \cite{GEZ2} and~\cite{EGZ2}, we
discussed the relations between these inverse principles and
similar principles formulated  by Arak
(see \cite{Arak0}--\cite{Arak and Zaitsev}) in his papers from the 1980's. In the one-dimensional
case, Arak has found a connection of the concentration function of
the sum with the arithmetic structure of supports of distributions
of independent random variables for {\it arbitrary}\/
distributions of summands. Using these results, he has solved an old problem
stated by Kolmogorov~\cite{K}.

In the present paper, we show that a consequence
of Arak's inequalities provides
results in the Littlewood--Offord problem of greater generality and improved precision compared to
those proved in \cite{GEZ2} and~\cite{EGZ2}. Moreover, using the
results of Tao and Vu \cite{TV08}, we
are able to describe the approximating sets much more precisely.

\medskip

Let us introduce first the necessary notations.  Below
${\mathbf N}$ and ${\mathbf N}_0$ {will denote} the sets of all
positive and non-negative integers respectively. The symbol $c$
will be used for absolute positive constants.
 Note that $c$ can be different in different (or even in the same) formulas.
We will write $A\ll B$ if $A\leq c B$. {Furthermore,  we
{will} use the notation} $A\asymp B$ if $A\ll B$ and $B\ll A$.
If the corresponding constant depends on, say, $r$, we write
$A\ll_r B$ and $A\asymp_r B$.  If $\xi=(\xi_1,\ldots,\xi_d)$ is a
random vector with distribution~$F=\mathcal L(\xi)$, we denote
$F^{(j)}=\mathcal L(\xi_j)$, $j=1,\ldots,d$. Let
$\widehat F(t)={\mathbf E}\,\exp\big(\,i\,\langle
t,\xi\rangle\big)$, \,$t\in\mathbf R^d$, \, be the characteristic
function of the distribution~$F$. Here $\langle\, \cdot \,,\,
\cdot \,\rangle$ is the inner product in \,$\mathbf R^d$.

For~${x=(x_1,\dots,x_n )\in\mathbf R^n}$ we denote
$\|x\|^2=\langle x,x\rangle= x_1^2+\dots +x_n^2$ and $|x|=
\max_j|x_j|$. Let \,$E_a$ \,be the distribution concentrated at a
point $a\in \mathbf{R}^n $. We denote by $[B]_\tau$ the closed
$\tau$-neighborhood of a set $B$ in the sense of the norm
$|\,\cdot\,|$. Products and powers of measures will be understood
in the sense of convolution. Thus, we write $F^n$ for the $n$-fold
convolution of a measure~$F$. While a distribution~$F$ is
infinitely divisible, $F^\lambda$, $\lambda\ge0$, is the
infinitely divisible distribution with characteristic function
$\widehat F^\lambda(t)$. For a finite set $K$, we denote by $|K|$
the number of elements~$x\in K$. \bl {The symbol $\times$} is used
for the direct product of sets. We write $O(\,\cdot\,)$ if the
involved constants depend on the parameters named ``constants'' in
the formulations, but not on~$n$.

The elementary properties of concentration functions are well
studied (see, for instance,
 \cite{Arak and Zaitsev, Hengartner and Theodorescu,
Petrov}). In particular, it is clear that
\begin{equation}\label{8j}Q(F,\mu)\le \big(1+\lfloor
\mu/\lambda\rfloor\big)^d\,Q(F,\lambda), \quad\hbox{for any }
\mu,\lambda>0,\end{equation} where $\lfloor x\rfloor$ is the
largest integer~$k$ that satisfies the inequality $k\le x$. Hence,
\begin{equation}\label{8a} Q(F,c\,\lambda)\asymp_d\,Q(F,\lambda).
\end{equation}
\medskip

Estimating the concentration functions in the Littlewood--Offord
problem, it is useful to reduce the problem to the estimation of
concentration functions of some symmetric infinitely divisible
distributions. The corresponding statement is contained in Lemma
\ref{lm42} below.

 Introduce the distribution
$H$ with the characteristic function
\begin{equation} \label{11}\widehat{H}(t)
=\exp\Big(-\frac{\,1\,}2\;\sum_{k=1}^{n}\big(1-\cos\left\langle \,
t,a_k\right\rangle \big)\Big).
\end{equation}
It is clear that $H$ is a symmetric infinitely divisible
distribution.
  Therefore, its characteristic function is positive for all $t\in \mathbf{R}^d$.

Let $\widetilde{X}=X_1-X_2$ \,be the symmetrized random vector,
where $X_1$ and $X_2$ are i.i.d.\ vectors involved in the
definition of $S_a$. In the sequel we use the notation
$G=\mathcal{L}(\widetilde{X})$. For $\delta\ge0$, we denote
\begin{equation}\label{pp}p(\delta)= G\big\{\{z:|z| >
\delta\}\big\}.\end{equation} Below we will use  the condition
\begin{equation}
G\{\{x\in\mathbf{R}:C_1<|x| < C_2\}\}\ge C_3, \label{1s8}
\end{equation}
where the values of $C_1, C_2, C_3$ will be specified in the
formulations below.

\begin{lemma}\label{lm42}
For any $\varkappa,\tau>0$, we have
\begin{equation}\label{1155}
Q(F_a, \tau) \ll_d Q(H^{p(\tau/\varkappa)}, \varkappa).
\end{equation}
\end{lemma}

According to \eqref{8j}, Lemma~\ref{lm42} implies the following
inequality.

\begin{corollary}\label{lm429}For any $\varkappa,\tau, \delta>0$, we have
\begin{equation}\label{1166}
Q(F_a, \tau) \ll_d\big(1+\lfloor\varkappa/\delta\rfloor \big)^d\,
Q(H^{p(\tau/\varkappa)}, \delta).
\end{equation}
\end{corollary}

Note that, in the case $\delta=\varkappa$, Corollary~\ref{lm42}
turns into Lemma~\ref{lm42}. Sometimes, it is useful to be free in
the choice of $\delta$ in~\eqref{1166}. In a recent paper of
Eliseeva and Zaitsev \cite{Eliseeva and Zaitsev2}, a more general
statement than Lemma \ref{lm42} is obtained. It gives useful
bounds if $p(\tau/\varkappa)$ is small, even if
$p(\tau/\varkappa)=0$. The proof of Lemma~\ref{lm42} is given
in~\cite{EGZ2}. It is rather elementary and is based on known
properties of concentration functions. We should note that
$H^{\lambda}$, $\lambda\ge0$, is a symmetric infinitely divisible
distribution with the L\'evy
 spectral measure $M_\lambda=\frac{\,\lambda\,}4\;M^*$, where
 $M^*=\sum_{k=1}^{n}\big(E_{a_k}+E_{-a_k}\big)$.

 Passing {to the limit $\tau \to 0$} in~\eqref{1155}, we obtain the following statement (see Zaitsev \cite{Zaitsev2}
 for details).

\begin{lemma}\label{lm342} The inequality
\begin{equation}\label{978t}
Q(F_a, 0) \ll_d Q(H^{p(0)}, 0) = H^{p(0)}\{\{0\}\}
\end{equation} holds.
\end{lemma}

Note that the case where $p(0)=0$ is trivial, since then
 $Q(F_a, 0) = Q(H^{p(0)}, 0) =1$ for any~$a$. Therefore we assume below that $p(0)>0$.

 \bigskip
The following definition is given in Tao and Vu \cite{TV08} (see
also \cite{Nguyen and Vu}, \cite{Nguyen and Vu13}, \cite{TV06},
and~\cite{Tao and Vu}).

Let $r\in{\mathbf N}_0$ be a non-negative integer,
$L=(L_1,\ldots, L_r)$
 be a $r$-tuple of positive reals,
and  $g=(g_1,\ldots, g_r)$ be a $r$-tuple of elements of $\mathbf
R^d$. The triplet $ P = (L, g, r)$  is called symmetric
{'generalized arithmetic progression' (GAP)} in $\mathbf R^d$. Here
$r$ is the rank, $L_1,\ldots, L_r$ are the dimensions and
$g_1,\ldots, g_r$ are the generators of the GAP~$P$. We define the
image $\hbox {Image\,}(P) \subset\mathbf R^d$ of $P$ to be the set
$$\hbox {Image}(P) = \big\{ m_1g_1 + \cdots + m_rg_r:-L_j \le m_j \le L_j,
 \ m_j\in{\mathbf Z} \hbox{ \ for all }j=1, \ldots,  r\big\}.
$$
  For $t>0$ we denote the dilate $P^t$ of $P$ as the symmetric
GAP $ P^t = (tL, g, r)$ with
$$\hbox {Image\,}(P^t) = \big\{ m_1g_1 + \cdots + m_rg_r:-tL_j \le m_j \le tL_j,
 \ m_j\in{\mathbf Z} \hbox{ \ for all }j=1, \ldots,  r\big\}.
$$
We define the size of $P$ to be $\hbox{size}(P) = \left|\,\hbox
{Image}(P)\right|$.

In fact, $\hbox {Image}(P)$ is the image of an integer box $$B =
\big\{(m_1,\ldots, m_r) \in{\mathbf Z}^r: -L_j \le m_j \le
L_j\big\}$$ under the linear map
$$
\Phi:(m_1,\ldots m_r) \in{\mathbf Z}^r \rightarrow  m_1g_1 +
\cdots + m_rg_r.
$$
We say that $P$ is proper if this map is one to one, or,
equivalently, if
\begin{equation}\label{6u2}\hbox{size}(P)=\prod_{j=1}^{r}\big(2\,\lfloor L_j\rfloor+1\big).
\end{equation}
The right-hand side of \eqref{6u2} is denoted by $\hbox{Vol}(P)$. It is called the volume of~$P$.
For non-proper GAPs, we  have, of course,
\begin{equation}\label{6u277}\hbox{size}(P)<\prod_{j=1}^{r}\big(2\,\lfloor
L_j\rfloor+1\big).
\end{equation}
For $t>0$, we say that $P$ is $t$-proper if $P^t$ is proper. It is
infinitely proper if it is $t$-proper for any $t>0$. In general,
for $t>0$, we have\begin{equation}\label{6u22}
\hbox{size}(P^t)\le\prod_{j=1}^{r}\big(2\,\lfloor
tL_j\rfloor+1\big).
\end{equation}

\begin{remark}\label{rem08}\rm In the  case $r=0$ the vectors $L$ and
$g$ have no elements and the image of the GAP~$P$ consists of the
unique zero vector $0\in\mathbf R^d$.
\end{remark}

\begin{remark}\label{rem00}\rm Symmetric GAPs are defined not only
by their images (sets of points in $\mathbf R^d$ admitting the
representation $m_1\,g_1 + \cdots + m_r\,g_r$, where $-L_j \le m_j \le
L_j$,
 $m_j\in{\mathbf Z}$, for $1 \le j\le  r$, see \cite{TV08}). The definition includes
 the generators $g_1,\ldots, g_r\in\mathbf R^d$ and the dimensions $L_1,\ldots,
L_r\in\mathbf R$. Different symmetric GAPs can have the same
image. For example, if $L_j<1$, then their generators $g_j$ are not
used in constructing the image of~$P$. However, the image of
~$P^t$ depends on~$g_j$ \;if \;$tL_j\ge1$. Obviously, by
definition, the GAPs $P$ and $P^t$ have the same generators and
the same rank.
\end{remark}

Recall that a convex body in the $r$-dimensional Euclidean space
$\mathbf R^r$ is a compact convex set with non-empty interior.

\begin{lemma}\label{lmj}
 Let $V$ be a convex symmetric body in
$\mathbf R^r$, and let $\Lambda$ be a lattice in $\mathbf R^r$.
Then there exists a symmetric, infinitely proper GAP $P$ in
$\Lambda$ with rank $l\le r$ such that we have
\begin{equation}
\label{6s2}\hbox {\rm Image}(P) \subset V\cap\Lambda\subset\hbox
{\rm Image} \bigl(P^{(c_1r)^{3r/2}}\bigr)
\end{equation}with an
absolute constant~$c_1\ge1$. Moreover, the generators $g_j$ of
$P$, for $1 \le j\le  l$, are contained in the symmetric body $l\hskip1pt V$.
\end{lemma}

The main part of Lemma \ref{lmj} is contained in Theorem~1.6 of
\cite{TV08}. The last statement of this Lemma follows from \cite[Theorem
3.34]{TV06}. The basis $g_1,\ldots, g_l\in\mathbf R^r$ is
sometimes called Mahler basis for the sublattice of $\Lambda$
spanned on  $\Lambda\cap V$.

\begin{corollary}\label{lmjj}Under  the conditions of Lemma\/
$\ref{lmj}$,
 \begin{equation} \hbox{\rm size}(P^{(c_1r)^{3r/2}})\le\big(2\,(c_1\,r)^{3r/2}+1\big)^r\left|\,V\cap\Lambda\,\right|.\label{7s2}
\end{equation}
\end{corollary}

\noindent {\it Proof of Corollary\/ $\ref{lmjj}$.} Using Lemma
\ref{lmj} and relations \eqref{6u2} and \eqref{6u22}, we obtain
\begin{eqnarray}
\hbox{size}\bigl(P^{(c_1r)^{3r/2}}\bigr)&\le&
\prod_{j=1}^{r}\big(\lfloor2\,
(c_1\,r)^{3r/2}L_j\rfloor+1\big)\nonumber\\
&\le&\big(2\,(c_1\,r)^{3r/2}+1\big)^r\hbox{size}(P)\le
\big(2\,(c_1\,r)^{3r/2}+1\big)^r\left|\,V\cap\Lambda\,\right|.\label{7y2}
\end{eqnarray}
Here the numbers $L_j$ are the dimensions of $P$. We used that
$\lfloor2\,tL\rfloor+1\le (2\,t+1)\,(\lfloor2\,L\rfloor+1)$, for
$L,t>0$. $\square$\medskip

The following Lemma~\ref{lmj4} shows that symmetric progressions
are contained in proper progressions. It can be found in
\cite[Theorem 1.9]{TV08}, see also \cite[Theorem 3.40]{TV06} and
\cite[Theorem~2.1]{Green}.
\begin{lemma}\label{lmj4}  Let $P$ be a symmetric GAP in
$\mathbf R$, and let $t \ge1$. Then there exists a $t$-proper
symmetric GAP $Q$ with $\hbox{\rm rank}(Q)\le r=\hbox{\rm rank}(P)$,
${\rm Image}(P) \subset {\rm Image}(Q)$, and
\begin{equation} \hbox{\rm size}(P)\le\hbox{\rm size}(Q)\le(2\,t)^{r}r^{6r^2}\hbox{\rm size}(P).\label{76s2}
\end{equation}
\end{lemma}\bigskip

We start now to formulate Theorem~\ref{thm2} which is a
one-dimensional Arak type result, see~\cite{Arak}. Let us
introduce the necessary notation.

Let  $r\in{\mathbf N}_0$, $m\in{\mathbf N}$ be fixed, let $h$
be an arbitrary $r$-dimensional vector,  and let $V$ be an
arbitrary closed symmetric convex subset of~${\mathbf R}^r$
containing not more than $m$ points with integer coordinates. We
define $\mathcal{K}_{r,m}$ as the collection of all sets of the
form
\begin{equation}
K=\big\{\langle{\nu}, h\rangle:{\nu}\in {\mathbf Z}^r\cap
V\big\}\subset{\mathbf R} .\label{1s1}
\end{equation}
We shall call such sets CGAPs ('convex generalized arithmetic
progressions', see \cite{Green}), by analogy with the notion of
GAPs.

Here the number~$r$ is the rank and $\left|\,{\mathbf Z}^r\cap
V\right|$ is the size of a CGAP  in the class $\mathcal{K}_{r,m}$.
It seems natural to call a CGAP from $\mathcal{K}_{r,m}$  {\it
proper} if all points $\big\{\langle{\nu}, h\rangle:{\nu}\in
{\mathbf Z}^r\big\}$ are disjoint.

For any Borel measure $W$ on ${\mathbf R}$ and $\tau\ge0$ we define $\beta_{r,m}(W, \tau)$ by
\begin{equation}
\beta_{r,m}(W, \tau)=\inf_{K\in\mathcal{K}_{r,m} }W\{{\mathbf
R}\setminus[K]_\tau\}. \label{1s3}
\end{equation}

We now introduce  a class  of $d$-dimensional CGAPs
$\mathcal{K}_{r,m}^{(d)}$ which consists of all sets of the form
$K=\times_{j=1}^{d} K_j$, where $K_j\in\mathcal{K}_{r_j,m_j}$,
$r=(r_1, \ldots,r_d)\in{\mathbf N}_0^d$, $m=(m_1,
\ldots,m_d)\in{\mathbf N}^d$. We call $R=r_1+\cdots +r_d$ the rank
and $\left|\,{\mathbf Z}^{r_1}\cap V_1\right|\cdot\; \cdots\;\cdot
\left|\,{\mathbf Z}^{r_d}\cap V_d\right|$ the size of~$K$. Here
$V_j\subset {\mathbf R}^{r_j}$ are symmetric convex subsets from
the representation \eqref{1s1} for $K_j$.

\begin{remark}\label{rem88}\rm In the  case $r=0$ the class $\mathcal{K}_{r,m}=\mathcal{K}_{0,m}$
consists of the one set $\{0\}$ having zero as the unique element.
\end{remark}

The following result is a particular case of Theorem~4.3 of Chapter~II  in \cite{Arak and Zaitsev}.

\begin{theorem}\label{thm2} Let $D$ be a one-dimensional
infinitely divisible distribution with characteristic function of
the form $\exp\big(\alpha\,(\widehat W(t)-1)\big)$, $t\in{\mathbf
R}$, where $\alpha>0$ and\/ $W$ is a probability distribution.
Let\/ $\tau\ge0$, \;$r\in{\mathbf N}_0$, $m\in{\mathbf N}$. Then
\begin{equation}
Q(D, \tau)\le
c_2^{r+1}\biggl(\frac{1}{m\sqrt{\alpha\,\beta_{r,m}(W, \tau)}}
+\frac{(r+1)^{5r/2}}{(\alpha\,\beta_{r,m}(W,
\tau))^{(r+1)/2}}\biggr), \label{1s4}
\end{equation}
where $c_2$ is an absolute constant.
\end{theorem}

\begin{remark}\label{rem0}\rm Arak \cite{Arak}
did not assume that the set $V$ is closed in the definition
\eqref{1s1}. It is easy to see however that this does not change
the formulation of Theorem~\ref{thm2},
\end{remark}

Arak \cite{Arak} proved an analogue of Theorem \ref{thm2} for sums
of i.i.d.\ random variables (see Theorem~4.2 of Chapter~II  in
\cite{Arak and Zaitsev}). He used this theorem in the proof of the
following remarkable result: \medskip

\noindent {\it There exists a universal constant\/ $C$ such that
for any one-dimensional probability distribution~$F$ and for any
positive integer~$n$ there exists an infinitely divisible
distribution~$D_n$ such that
$$
\rho(F^{n}, D_n)\le C\,n^{-2/3},
$$
where $ \rho (\, \cdot \,, \, \cdot \,) $ is the classical Kolmogorov's uniform distance between
corresponding distribution functions.
}\medskip

\noindent This gives the definitive solution of an old problem
stated by Kolmogorov~\cite{K} in the 1950's (see~\cite{Arak and
Zaitsev} for the history of this problem).

Estimation of concentration functions is the main tool for the bound of $\rho(F^{n}, D_n)$.
Moreover, Arak's inequalities play a crucial role in {\it constructing} the approximating distribution~$D_n$.
Standard distributions, such as Gaussian, stable, accompanying compound Poisson laws,
won't suffice to ensure the bound for arbitrary~$F$ and~$n$.
Roughly speaking, the main idea is that either the concentration function of $F^{n}$
is relatively small (and the standard approximation is good enough)
or the support of the distribution of the summands
has a simple arithmetical structure. In the latter case, this structure is used for
constructing the distribution~$D_n$.

The investigations of  Arak in~\cite{Arak0} and~\cite{Arak} were motivated  by the ideas
 of Freiman~\cite{16} on the structural theory of set addition. These ideas were used by
 Nguyen and Vu~\cite{Nguyen and Vu} and~\cite{Nguyen and Vu13} as well.
 The proof of Theorem~\ref{thm2} is based on Ess\'een's inequality \cite{14}
 estimating the concentration function by an integral of the modulus of characteristic function.
 The important tools used are the Parseval equality and
 the following obvious inequality for characteristic functions of  one-dimensional distributions~$U$:
$$
\bigl|\widehat U(t+h)-\widehat U(t)\bigr|^2\le 2\,\big(1-{\rm Re}\,\widehat U(h)),
\quad\hbox{for all }t,h\in \mathbf{R}.
$$
\medskip

Corollary~\ref{lm429}, Lemma~\ref{lm342} and
Theorem~\ref{thm2} imply the following Theorem~\ref{thm7}.

 \begin{theorem}\label{thm7} Let $\varkappa, \delta>0$,  $\tau\ge0$, and let $X$ be a real random variable satisfying
 condition~\eqref{1s8} with $C_1=\tau/\varkappa$, $C_2=\infty$ and $C_3=p(\tau/\varkappa)>0$.
Let $d=1$, \;$r\in{\mathbf N}_0$, $m\in{\mathbf N}$. Then
\begin{equation}
Q(F_a, \tau)\le c_3^{r+1}\big(1+\lfloor\varkappa/\delta\rfloor
\big)\,\biggl(\frac{1}{m\sqrt{p(\tau/\varkappa)\,\beta_{r,m}(M^*, \delta)}}
+\frac{(r+1)^{5r/2}}{(p(\tau/\varkappa)\,\beta_{r,m}(M^*, \delta))^{(r+1)/2}}\biggr),
\quad\hbox{if }\tau>0, \label{1sy4}
\end{equation}
and
\begin{equation}
Q(F_a, 0)\le c_3^{r+1}\,\biggl(\frac{1}{m\sqrt{p(0)\,\beta_{r,m}(M^*,
0)}} +\frac{(r+1)^{5r/2}}{(p(0)\,\beta_{r,m}(M^*, 0))^{(r+1)/2}}\biggr),
\quad\hbox{if }\tau=0, \label{1sy49}
\end{equation}
where
 $M^*=\sum_{k=1}^{n}\big(E_{a_k}+E_{-a_k}\big)$ and $c_3$ is an absolute constant.
\end{theorem}

In order to prove Theorem~{\ref{thm7}}, it suffices to apply
Corollary~\ref{lm429}, Lemma~\ref{lm342} and Theorem~\ref{thm2}
and to note that $H^{p(\tau/\varkappa)}$ is an infinitely
divisible distribution whose L\'evy
 spectral measure is $p(\tau/\varkappa)\,M^*/4$.
 Introduce as well $M=\sum_{k=1}^{n}E_{a_k}$.
 It is obvious that $M\le M^*$ and $\beta_{r,m}(M, \delta)\le\beta_{r,m}(M^*, \delta)$.
\medskip

\section{Results}\label{2}

The main results of the present paper are Theorems \ref{thm17} and
\ref{thm27}. Their proofs are based on Theorem \ref{thm7}.
Theorems \ref{thm7}, \ref{thm17} and \ref{thm27} have
non-asymptotic character.  They provide information about the
arithmetic structure of $a=(a_1,\ldots,a_n)$ without assumptions
like $q_j=Q(F_a^{(j)}, \tau)\ge n^{-A}$, $j=1,\ldots,d$, which are
imposed in Theorem \ref{nthm8} below. Theorems \ref{thm17} and
\ref{thm27} are formulated for {\it fixed~$n$} and the dependence
of constants on parameters is given explicitly. No analogues of Theorems \ref{thm17} and
\ref{thm27} follow from the asymptotical results of Nguyen, Tao and Vu \cite{Nguyen and Vu},
\cite{Nguyen and Vu13}
and~\cite{Tao and Vu}, see Theorem~\ref{tnv2}.  The conditions of Theorem
\ref{nthm8} are weaker than those used in the results
of Nguyen, Tao and Vu. Theorem \ref{nthm8}
was derived from Theorem~\ref{thm7} in the
 paper of Eliseeva, G\"otze and  Zaitsev \cite{EGZ2}.

In the following we state Theorems \ref{thm16} and \ref{thm19} which
are more general than Theorem \ref{nthm8}. Theorems \ref{thm16}
and~\ref{thm19} will be deduced from Theorem~\ref{thm17}. We
conclude with a comparison of Theorems \ref{nthm8}, \ref{thm16} and~\ref{thm19}
with the results of Nguyen, Tao and~Vu. Notice that in the
asymptotic Theorems~\ref{nthm8}--\ref{tnv2}, where $n\to\infty$,
the elements $a_{j}$ of\/ $a$ may depend on $n$.

Finally, we state improved and generalized versions of Theorems~5 and~6 of~\cite{EGZ2},
see Theorems~\ref{tnv3}--\ref{thm55} of the present paper.

 \begin{theorem}\label{thm17}
Let\/  $d=1$,  \;$a=(a_1,\ldots,a_n) \in {\mathbf{R}^n}$, \;$p(0)>0$,
and\/ $q=Q(F_a, \tau)$,   $\tau\ge0$. There exists positive
absolute constants $c_4$--$c_7$
 such that for any\/ $\varkappa>0$,  $\delta\ge0$, for any fixed $r\in{\mathbf N}_0$,
 and any $n'\in \mathbf N$ satisfying the inequalities $\delta\le\max\{\varkappa, \tau\}$ and
\begin{equation}\big(\,2\,c_4^{r+1}\,(r+1)^{5r/2}\,\varkappa\big/q\,\delta\,\big)^{2/(r+1)}/p(\tau/\varkappa)
\le n' \le n, \quad\hbox{if }\tau>0,\label{2b1}
\end{equation}
or
\begin{equation}\big(\,2\,c_4^{r+1}\,(r+1)^{5r/2}\big/q\,\big)^{2/(r+1)}/p(0)
 \le n' \le n, \quad\hbox{if }\tau=0,\label{2b13}\end{equation}
 there exist $m\in \mathbf N$ and CGAPs $K^*, K^{**}\subset\mathbf R$ having ranks
 $\le r$ and sizes~$\ll m$ and $\ll c(r)\,m$ respectively and such that

$1$. At least $n-2\,n'$ elements $a_{k}$ of\/ $a$ are
$\delta$-close to $K^*$, that is, $a_{k}\in[K^*]_{\delta}$ $($this
means that for these elements $a_{k}$ there exist $y_{k}\in K^*$
such that $\left|a_{k}-y_{k}\right|\le\delta)$.

$2$. The above number $m$ satisfies the inequalities
 \begin{equation}
 m\le\frac{2\,c_4^{r+1}\,\varkappa}{ q\,\delta\,\sqrt{p(\tau/\varkappa)\,
n'}}+1, \quad\hbox{if }\tau>0, \label{2s4478}
\end{equation}
or
\begin{equation}
 m\le\frac{2\,c_4^{r+1}}{ q\,\sqrt{p(0)\,
n'}}+1, \quad\hbox{if }\tau=0. \label{2s447}
\end{equation}

$3$.  The set $K^*$ is contained in the image $\overline K$ of a
symmetric GAP $\overline P$ which has rank \,$\overline l\le r$,
 size~$\ll (c_5\,r)^{3r^2/2}m$ and generators $\overline g_j$,
$j=1,\ldots,\overline l$, satisfying inequality
 $\bigl|\overline g_j\bigr|\le 2\,r
\left\|a\right\|/\sqrt{n'}$.

$4$.  The set $K^*$ is contained in the image $\overline{\overline
K}$ of a proper symmetric GAP $\overline{\overline P}$ which has
rank \,$\overline{\overline l}\le r$ and size~$\ll
(c_6\,r)^{15r^2/2}m$.

$5$. At least $n-2\,n'$ elements of\/ $a$ are  $\delta$-close to
 $K^{**}$.

$6$.  The set $K^{**}$ is contained in the image $\widetilde K$ of
a proper symmetric GAP ${\widetilde P}$ which has rank
\,$\widetilde l\le r$, size~$\ll (c_7\,r)^{21r^2/2}m$ and
generators $\widetilde g_j$, $j=1,\ldots,\widetilde l$, satisfying
the inequality
 $\bigl|\widetilde g_j\bigr|\le 2\,r
\left\|a\right\|/\sqrt{n'}$.
\end{theorem}

The statement of Theorem~\ref{thm17} is rather cumbersome, but
this is the price for its generality. The formulation may be
simplified in particular cases, for example, for
$\varkappa=\delta$ or for $\varkappa=\tau$.

The assertion of Theorem~\ref{thm17} is
non-trivial for each fixed $r$ starting with $r=0$. In this
case $m=1$ and Theorem~\ref{thm17} gives a bound for the amount
$N$ of elements $a_{k}$ of~$a$ which are outside of the interval
$[-\delta,\,\delta]$ around zero.
Namely,
$$N\leq\big(\,2\,c_4\,\varkappa\big/q\,\delta\,\big)^{2}/p(\tau/\varkappa)+1,
\quad\hbox{if }\tau>0,$$ and
$$N\leq\big(\,2\,c_4\big/q\,\big)^{2}/p(0)+1, \quad\hbox{if }\delta=\tau=0.$$

Comparing item 3 with items 4 and 6 of Theorem~\ref{thm17}, we see
that in item 3 the approximating GAP may be non-proper. However,
the size of proper approximating GAPs is larger in items 4 and 6.
Moreover, if $\delta>0$, then it is obvious that by small
perturbations of generators of a non-proper GAP~$\overline P$ with
$\hbox{Image}(\overline P)=\overline K$, we can construct a proper
GAP~$\overline{\overline P}$ with
$\hbox{Image}(\overline{\overline P})=\overline{\overline K} $,
with the size and generators satisfying the bounds of item 3 and
 such that $[\overline K]_\delta\subset
[\overline{\overline K} ]_{2\delta}$. The set $
[\overline{\overline K}]_{2\delta}$ approximates the set of elements of~$a$ not
worse than $[{\overline K}]_\delta$. Note that, according to
\eqref{8a}, in the conditions of our results there is no essential
difference between $\delta$ and $2\delta$-neighborhoods. Thus, in
fact the statements of items 4--6 are useful in the case
$\delta=\tau=0$ only. Otherwise, the statement of item 3 is good
enough.

\begin{remark}\label{rem3}\rm Notice that Theorem~\ref{thm17} does not provide any information
if there is no $n'$ satisfying
inequalities~\eqref{2b1} or~\eqref{2b13}. In particular,
if $\big(\,2\,c_4^{r+1}\,(r+1)^{5r/2}\,\varkappa\big/q\,\delta\,\big)^{2/(r+1)}/p(\tau/\varkappa)
 \ge n$ and $\tau>0$. The same may be said if $n-2\,n'\leq0$.
 Similar remarks can be made about Theorems~\ref{thm27}--\ref{thm55}.
\end{remark}

Theorem~\ref{thm17} is  formulated for one-dimensional $a_k$, $k=1,\ldots, n$.
  However, it may be shown that
Theorem~\ref{thm17} provides sufficiently rich  arithmetic
properties for the set $a=(a_1,\ldots,a_n) \in {(\mathbf{R}^d)}^n$
in the multivariate case as well (see Theorem~\ref{thm27} below).
It suffices to apply Theorem~\ref{thm17} to the distributions
$F_a^{(j)}$, $j=1,\ldots,d$, where $F_a^{(j)}$ are distributions
of coordinates of the vector~$S_a$.\medskip

 Introduce the vectors $a^{(j)}=(a_{1j},\ldots,a_{nj})$,
$j=1,\ldots, d$. It is obvious that $F_a^{(j)}=F_{a^{(j)}}$.

\begin{theorem}\label{thm27}
Let\/  $d>1$,  \;$p(0)>0$, $q_j=Q(F_a^{(j)}, \tau_j)$,
$\tau_j\ge0$, $j=1,\ldots,d$. Below\/ $c_4$--$c_7$ are positive
absolute constants from Theorem~$\ref{thm17}$. Suppose that\/
$a=(a_1,\ldots,a_n) \in ({\mathbf{R}^d})^n$ is a multi-subset of\/
${\mathbf R}^d$.
 Let\/ $\varkappa_j>0$, $\delta_j\geq0$,  $r_j\in \mathbf N_0$,
 and\/  $n'_j\in \mathbf N$, $j=1,\ldots,d$, satisfy inequalities\/  $\delta_j\le\max\{\varkappa_j, \tau_j\}$ and
 \begin{equation}\big(\,2\,c_4^{r_j+1}\,(r_j+1)^{5r_j/2}\,
 \varkappa_j\big/q_j\,\delta_j\,\big)^{2/(r_j+1)}\big/p(\tau_j/\varkappa_j)\le n_j' \le n,
 \quad\hbox{if }\tau_j>0,\label{2s45}
\end{equation}
 or
\begin{equation}\big(\,2\,c_4^{r_j+1}\,(r_j+1)^{5r_j/2}\big/q_j\,\big)^{2/(r_j+1)}\big/p(0)\le n_j' \le n,
 \quad\hbox{if }\tau_j=0.\label{2s46}
\end{equation}
 Then, for each $j=1,\ldots,d$, there exist $m_j\in \mathbf
N$  and CGAPs $K_j^*, K_j^{**} \subset\mathbf R$ having ranks $\le
r_j$ and sizes~$\ll m_j$ and $\ll c(r_j)\,m_j$ respectively and
such that

$1$. At least $n-2\,n'_j$ elements $a_{kj}$ of\/ $a^{(j)}$ are
$\delta_j$-close to
 $K_j^*$, that is,
$a_{kj}\in[K_j^*]_{\delta_j}$
 $($this means that for these elements $a_{kj}$  there exist $y_{kj}\in K_j^*$ such
that $\left|a_{kj}-y_{kj}\right|\le\delta_j)$.

$2$. $m_j$ satisfies inequality $m_j\le w_j$, where
 \begin{equation}
 w_j=\frac{2\,c_4^{r_j+1}\,\varkappa_j}{ q_j\,\delta_j\,\sqrt{p(\tau_j/\varkappa_j)\,
n'_j}}+1,\quad\hbox{if }\tau_j>0, \label{2s498}
\end{equation}
 or
  \begin{equation}
 w_j=\frac{2\,c_4^{r_j+1}}{ q_j\,\sqrt{p(0)\,
n'_j}}+1,\quad\hbox{if }\tau_j=0. \label{2d8}
\end{equation}

$3$.  The set $K_j^*$ is contained in the image $\overline K_j$ of
a symmetric GAP $\overline P_j$ which has rank \,$\overline l_j\le
r_j$, size~$\ll (c_5\,r_j)^{3r_j^2/2}m_j$ and generators
$\overline g_p^{(j)}$, $p=1,\ldots,\overline l_j$, satisfying
inequality
 $\bigl|\overline g_p^{(j)}\bigr|\le 2\,r_j
 \left\|a^{(j)}\right\|/\sqrt{n'_j}$.

$4$.  The set $K_j^*$ is contained in the image
$\overline{\overline K}_j$ of a proper symmetric GAP
$\overline{\overline P}_j$ which has rank \,$\overline{\overline
l}_j\le r_j$ and size~$\ll (c_6\,r_j)^{15r_j^2/2}m_j$.

$5$.  At least $n-2\,n'_j$ elements of\/ $a^{(j)}$ are
$\delta_j$-close to
 $K_j^{**}$.

$6$.  The set $K_j^{**}$ is contained in the image $\widetilde
K_j$ of a proper symmetric GAP ${\widetilde P}_j$ which has rank
\,$\widetilde l_j\le r_j$, size~$\ll (c_7\,r_j)^{21r_j^2/2}m_j$
and generators $\widetilde g_p^{(j)}$, $p=1,\ldots,\widetilde
l_j$, satisfying inequality
 $\bigl|\widetilde g_p^{(j)}\bigr|\le 2\,r_j
\left\|a^{(j)}\right\|/\sqrt{n'_j}$.

\medskip

$7$. The multi-vector\/ $a$ is well approximated by
 the $d$-dimensional CGAPs $K^*=\times _{j=1}^{d} K_j^*$ and $K^{**}=\times _{j=1}^{d} K_j^{**}$, by the
image $\overline K=\times _{j=1}^{d}\overline K_j$ of a symmetric
GAP $\overline P$, and by the the images $\overline{\overline
K}=\times _{j=1}^{d} \overline{\overline
 K}_j$ and
 $\widetilde K=\times _{j=1}^{d}\widetilde  K_j$  of proper
 symmetric GAPs
 $\overline{\overline P}$ and
 $\widetilde P$ of ranks
 $\le R=r_1+\cdots
 +r_d$.
 At least $n-2\sum_{j=1}^dn'_j$ elements of\/ $a$ \,belong to each
 of the sets $\times_{j=1}^{d}[K^*_j]_{\delta_j}$,
 $\times_{j=1}^{d}[K^{**}_j]_{\delta_j}$, $\times_{j=1}^{d}[{\overline K}_j]_{\delta_j}$,
 $\times_{j=1}^{d}[\overline{\overline K}_j]_{\delta_j}$,
 and\/ $\times_{j=1}^{d}[{\widetilde K}_j]_{\delta_j}$.
Furthermore,
\begin{equation}\label{32dd}\left|\times _{j=1}^{d} K_j^*\right|\ll_d \prod_{j=1}^{d} m_j
\le\prod_{j=1}^{d}
 w_j ,\end{equation}
\begin{equation}\label{32ed}\left|\times _{j=1}^{d} {\overline K}_j\right|\ll_d \prod_{j=1}^{d} (c_5\,r_j)^{3r_j^2/2}\,m_j
\le\prod_{j=1}^{d}(c_5\,r_j)^{3r_j^2/2}\,
 w_j ,\end{equation}
\begin{equation}\label{32fd}\left|\times _{j=1}^{d} \overline{\overline K}_j\right|
\ll_d \prod_{j=1}^{d}(c_6\,r_j)^{15r_j^2/2}\, m_j
\le\prod_{j=1}^{d}(c_6\,r_j)^{15r_j^2/2}\,
 w_j ,\end{equation}and
 \begin{equation}\label{132ed4}\left|\times _{j=1}^{d} {\widetilde K}_j\right|\ll_d \prod_{j=1}^{d} (c_7\,r_j)^{21r_j^2/2}\,m_j
\le\prod_{j=1}^{d}(c_7\,r_j)^{21r_j^2/2}\,
 w_j ,\end{equation}

$8$.  Let\/ $l_j$ be a short notation for\/ ${\overline
l}_j,\overline{\overline l}_j,{\widetilde l}_j$, the 
{number} of generators of ${\overline P}_j,\overline{\overline
P}_j,{\widetilde P}_j$ respectively. The generators ${\overline
g}_s,\overline{\overline g}_s,{\widetilde g}_s\in {\mathbf R}^d$,
$s=1,\ldots,l_1+\cdots+l_d$, of the GAPs~${\overline P}$,
$\overline{\overline P}$, ${\widetilde P}$ respectively have only
one non-zero coordinate each. Denote
$$
 s_0=0\quad\hbox{and}\quad s_k=\sum_{j=1}^{k}l_j, \quad k=1,\ldots,d.
$$
For $s_{k-1}<s\le s_k$, the generators ${\overline
g}_s,\overline{\overline g}_s,{\widetilde g}_s$ are non-zero in
the $k$-th coordinates only and these coordinates are equal to the
sequence of  generators ${\overline g}_1^{(k)},\ldots, {\overline
g}_{l_k}^{(k)}$; $\overline{\overline g}_1^{(k)},\ldots,
\overline{\overline g}_{l_k}^{(k)}$; ${\widetilde
g}_1^{(k)},\ldots, {\widetilde g}_{l_k}^{(k)}$ of the
GAPs~${\overline P}_k$, $\overline{\overline P}_k$, ${\widetilde
P}_k$ respectively satisfying inequality
$$\max\bigl\{\bigl|\overline g_p^{(k)}\bigr|, \bigl|\widetilde
g_p^{(k)}\bigr| \bigr\}\le 2\,r_k
 \left\|a^{(k)}\right\|/\sqrt{n'_k},
\quad p=1,\ldots,l_k.$$
\end{theorem}

\medskip

The following Theorem \ref{nthm8} was obtained in  \cite{EGZ2}
with the use of Theorem~\ref{thm7} (see \cite[Theorem~3 and
Proposition~1]{EGZ2}).

 \begin{theorem}\label{nthm8}
Let $d\ge1$,  $0 < \varepsilon \le 1$, $0 < \theta \le 1$, $A>0$,
$B>0$, $C_3>0$ be constants and $\tau_{n} \geq 0$ be a parameter
that may depend on $n$.  Let $X$ be a real random variable
satisfying
 condition~\eqref{1s8} with $C_1=1$, $C_2=\infty$ and $C_3\le p(1)$.
 Suppose that $a=(a_1,\ldots,a_n) \in
{(\mathbf{R}^d)}^n$ is a multi-subset of\/~${\mathbf R}^d$ such
that $q_j=Q(F_a^{(j)}, \tau_n)\ge n^{-A}$, $j=1,\ldots,d$, where
$F_a^{(j)}$ are distributions of coordinates of the vector~$S_a$.
 Let \,$\rho_n$
denote  a non-random sequence satisfying $ n^{-B}\leq\rho_n\leq1$
and let \,$\delta_{n}=\tau_{n}\,\rho_n$.
 Then, for any number $n'$ such that \,$\varepsilon \,n^\theta\le n'\leq n$,
 there exists a proper symmetric GAP $P$ such that

 $1$. At least $n-dn'$ elements $a_{j}$ of\/ $a$ are
$\delta_n$-close to the image $K$ of the GAP $P$ in the norm
$|\,\cdot\,|$.

 $2$. $P$ has small rank $R=O(1)$, and small size
 \begin{equation}
\hbox{\rm size}(P)=|K|\le
\prod_{j=1}^d\max\Big\{O\Big(q_j^{-1}\,\rho_n^{-1}\,(n')^{-1/2}\Big),
1\Big\}. \label{n11sp}
\end{equation}
\end{theorem}

\begin{remark}\label{rem}\rm In the first version of the preprint of the paper~\cite{EGZ2},
the GAP $K$ may be non-proper in Theorem~\ref{nthm8}. In order to
get the properness of $K$ we have
moreover used
arguments from a paper of Tao and Vu \cite{TV08}.
\end{remark}

Theorems \ref{thm16} and~\ref{thm19} below are consequences of
Theorem~\ref{thm17}.

 \begin{theorem}\label{thm16}
 Let\/ $b_n>0$, $n=1, 2, \ldots,$ \,be a $($depending on $n)$ sequence of non-random parameters
tending to infinity as $n\to\infty$. Let  $A, \theta,
\varepsilon_1,\varepsilon_2 >0$ be constants,  and\/ \;$p(0)>0$.
Suppose that $a=(a_1,\ldots,a_n) \in {(\mathbf{R}^d)}^n$ is a
multi-subset of\/~${\mathbf R}^d$ such that $q_j=Q(F_a^{(j)},
0)\ge \varepsilon_1\,b_n^{-A}$, $j=1,\ldots,d$, where $F_a^{(j)}$
are distributions of coordinates of the vector~$S_a$.
 Then, for any numbers\/ $n'_j$ such that\/ $\varepsilon_2\,b_n^\theta\leq n'_j\leq n$, $j=1,\ldots,d$,
 there exists a proper symmetric GAP\/ $P$ such that\medskip

$1$. At least $n-2\sum_{j=1}^dn'_j$ elements of\/ $a$ belong to
the image $K$ of the GAP $P$.
\medskip

 $2$. $P$ has small rank $L\le R=r_1+\cdots +r_d=O(1)$, and small size
 \begin{equation}
\hbox{\rm size}(P)=|K|\le
\prod_{j=1}^d\max\Big\{O\Big(q_j^{-1}\,(n'_j)^{-1/2}\Big),
1\Big\}. \label{1n11sp5}
\end{equation}\medskip

 $3$. Moreover, $K=\times _{j=1}^{d} K_j$, where $K_j$ are images
 of one-dimensional symmetric GAPs $P_j$ of rank $l_j\le r_j$, and the generators $g_s$, $s=1,\ldots,L=l_1+\cdots +l_d$,
of the GAP~$P$ of rank $L$ have only one non-zero coordinate each.
Denote
$$
 s_0=0\quad\hbox{and}\quad s_k=\sum_{j=1}^{k}l_j, \quad k=1,\ldots,d.
$$
For $s_{k-1}<s\le s_k$, the generators ${g}_s$ are non-zero in the
$k$-th coordinates only and these coordinates are equal to the
sequence of  generators ${g}_1^{(k)},\ldots, {g}_{l_k}^{(k)}$ of
the GAPs~$P_k$,  satisfying the inequality
 $\bigl|g_p^{(k)}\bigr|\le
2\,r_k\, \left\|a^{(k)}\right\|/\sqrt{n'_k}$, \;$p=1,\ldots,l_k$.
\end{theorem}\medskip

 \begin{theorem}\label{thm19}
 Let\/  $\theta, A,\varepsilon_1,\varepsilon_2,\varepsilon_3,\varepsilon_4>0$ and\/ $B,D\ge0$
be constants, and\/ $\theta>D$. Let\/ $b_n,\varkappa_n, \delta_n,
\tau_n,\rho_n>0$, $n=1,2,\ldots$, be depending on $n$ non-random
parameters satisfying the relations \;$p(\tau_n/\varkappa_n)\ge
\varepsilon_3\,b_n^{-D}$, \;$
\varepsilon_4\,b_n^{-B}\leq\rho_n=\delta_n/\varkappa_n\leq1$,
$\delta_n\le\max\{\varkappa_n, \tau_n\}$, for all $n\in \mathbf
N$, and\/ $b_n\to\infty$ as $n\to\infty$.
 Suppose that $a=(a_1,\ldots,a_n) \in
{(\mathbf{R}^d)}^n$ is a multi-subset of\/~${\mathbf R}^d$ such
that $q_j=Q(F_a^{(j)}, \tau_n)\ge \varepsilon_1\,b_n^{-A}$,
\;$j=1,\ldots,d$, \;$n=1,2,\ldots$, where $F_a^{(j)}$ are
distributions of coordinates of the vector~$S_a$.
 Then, for any numbers $n'_j$ such that $\varepsilon_2\,b_n^\theta\leq n'_j\leq n$, $j=1,\ldots,d$,
 there exists a  proper symmetric GAP $P$ such that

$1$. At least $n-2\sum_{j=1}^dn'_j$ elements of\/ $a$ are
$\delta_n$-close to the image $K$ of the GAP $P$ in the norm
$|\,\cdot\,|$.

 $2$. $P$ has small rank $L\le R=r_1+\cdots +r_d=O(1)$, and small size
 \begin{equation}
|K|\le
\prod_{j=1}^d\max\Big\{O\Big(q_j^{-1}\,\rho_n^{-1}\,(n'_j\,p(\tau_n/\varkappa_n))^{-1/2}\Big),
1\Big\}. \label{n111sp}
\end{equation}\medskip

 $3$.  Moreover, the properties of the
generators $g_s$ of the GAP~$P$ described in the item~$3$ of the
formulation of Theorem\/ $\ref{thm16}$ are still satisfied. In
particular, the inequality
 $\bigl\|g_s\bigr\|\le
2\,r_k\, \left\|a^{(k)}\right\|/\sqrt{n'_k}$ hold, for
$s_{k-1}<s\le s_k$.
\end{theorem}
\medskip

Theorem~\ref{thm19} is more general than Theorem~\ref{nthm8},
where we restricted ourselves to the case $b_n=n$, $\varkappa_n=\tau_n$, {$n_j'=n'$} only.
Theorem~\ref{thm19}
 provides  new substantial information if, for instance, the
 ratio $\tau_n/\delta_n$ is large and if $p(\tau_n/\varkappa_n)$ is not too small.

In applications of Theorem~\ref{thm19}, it is sometimes useful to minimize
the parameter $\delta_n$ responsible for the size of the
neighborhood of the set~$K$. Assume, for simplicity, that
$\delta_n=\varkappa_n$, for all $n\in \mathbf N$. Then the
condition \;$p(\tau_n/\delta_n)=p(\tau_n/\varkappa_n)\ge
\varepsilon_3\,b_n^{-D}$ is satisfied for larger values of
$\tau_n/\delta_n$ if the function $p(x)$ decreases slowly as $x\to
\infty$, that is, if the distribution $\mathcal L(\widetilde X)$
has heavy tails. Moreover, it is clear that for any function
$f(n)$ tending to infinity as $n\to \infty$ there exists a
distribution $\mathcal L(\widetilde X)$ such that
\;$p(\tau_n/\delta_n)\ge \varepsilon_3\,b_n^{-D}$ and
\,$\tau_n/\delta_n\ge f(n)$, for sufficiently large~$n$.\bigskip

A discussion concerning the comparison of Theorem~\ref{nthm8} with
the results of Nguyen, Tao and Vu \cite{Nguyen and Vu, Nguyen and
Vu13, Tao and Vu, Tao and Vu3} is given in \cite{GEZ2} and
\cite{EGZ2}.

 A few years ago Tao and Vu
\cite{Tao and Vu} formulated in the discrete case (with
$\tau_{n}=0)$ the so-called {'inverse principle'}, stating that
$$\hbox{\it A set $a=(a_1,\ldots,a_n)$ with large small
ball probability must have strong additive structure.}$$ Here
``large small ball probability'' means that
$Q(F_a,0)=\max_x{\mathbf P}\{S_a=x\}\ge n^{-A}$ with some constant
$A>0$. ``Strong additive structure'' means that a large part of
vectors $a_1,\ldots,a_n$ belong to a GAP with bounded size.

 Nguyen and Vu \cite{Nguyen and Vu}
have extended {this inverse principle} to the continuous case (with
$\tau_{n}>0)$ proving, in particular, the following result.

\begin{theorem}\label{tnv2}
Let $X$ be a real random variable satisfying condition \eqref{1s8}
with positive constants $C_1,C_2,C_3$. Let $0 < \varepsilon < 1$,
$A>0$ be constants and $\tau_{n} > 0$ be a parameter that may
depend on $n$. Suppose that $a=(a_1,\ldots,a_n) \in
({\mathbf{R}^d})^n$ is a multi-subset of ${\mathbf R}^d$ such that
  $q=Q(F_a, \tau_{n})\ge n^{-A}$.
 Then,
for any number $n'$ between $n^\varepsilon$ and $n$,
 there exists a symmetric proper GAP $P$ with image
 $K$ such that

 $1$. At least $n-n'$ elements of $a$ are  $\tau_{n}$-close to $K$.

 $2$. $P$ has small rank $r=O(1)$, and small size
 \begin{equation}
|K|\le \max\big\{O\big(q^{-1}(n')^{-1/2}\big), 1\big\}.
\label{12sp}
\end{equation}

$3$. There is a non-zero integer $p = O(\sqrt{n'})$
 such that all generators $g_j$ of $P$ have the form
$g_j = (g_{j1},\ldots,g_{jd})$, where $g_{jk} = \|a\|\,\tau_{n}\,
p_{jk}/p$
 with  \;$p_{jk}\in\mathbf Z$ and $p_{jk} = O(\tau_{n}^{-1}\sqrt{n'})$.
 \end{theorem}

Theorem~\ref{nthm8}
 allows us to derive  a
one-dimensional version of the first two statements of
Theorem~\ref{tnv2}. Theorem~\ref{thm19} contains an analogue of
the third one. Moreover, in Theorem~\ref{thm19}, the generators of
approximating GAPs have norms bounded from above by the quantities
with $\sqrt{n_{k}'}$ in the denominator, in contrast with Theorem
\ref{tnv2}. Furthermore, the condition $Q(F_a, \tau_n)\ge n^{-A}$
of Theorem~\ref{tnv2} implies the condition $Q(F_a^{(j)},
\tau_n)\ge n^{-A}$, $j=1,\ldots,d$, of Theorem~\ref{nthm8}, since
$Q(F_a^{(j)}, \tau_n)\ge Q(F_a, \tau_n)$. In addition,
$C_1,C_2,C_3$ are finite fixed constants in Theorem~\ref{tnv2},
while $C_2=\infty$ in Theorems~\ref{thm7}--\ref{thm19},
 and $C_1=\tau_n/\varkappa_n$, $C_3= \varepsilon_3\,b_n^{-D}$ in Theorem~\ref{thm19}. Notice
 that $p(1)$ and $p(\tau_{n}/\varkappa_{n})$ are involved in our Theorems~\ref{nthm8} and~\ref{thm19}
 explicitly. Theorem~\ref{tnv2} corresponds to the case $b_n=n$ in the
more general Theorems~\ref{thm16} and~\ref{thm19}, in which, however,
 number-theoretical properties of generators (as in Theorem~\ref{tnv2}, item~3) are not provided .
 We would like to emphasize the
non-asymptotic character of Theorems \ref{thm7}, \ref{thm17} and
\ref{thm27}.

For technical reasons,  we use in our results the quantity
$2\,n'$ instead of $n'$ for the number of points which can be not
approximated. It is clear that this difference is not significant.

We have to say that there
are some results from \cite{Nguyen and Vu, Nguyen and Vu13, Tao
and Vu, Tao and Vu3} which do not follow from the results of Arak.
In particular, we don't consider distributions on general additive
groups.

Sometimes, for $d>1$, inequality \eqref{n11sp} (with
$\rho_{n}=1$) or inequality~\eqref{n111sp} (with
$\rho_{n}=1$, $\varkappa_{n}=\tau_{n}$, $b_{n}={n}$) may be even stronger than inequality~\eqref{12sp}.
For example, if the vector $S_a$ has independent coordinates (this
may happen if each of the vectors~$a_j$ has only one non-zero
coordinate), then
 \begin{equation}
q=Q(F_a,\tau)\asymp_d \prod_{j=1}^d q_j.
 \label{11sps}
\end{equation}
Assuming, for simplicity, that $q_1=\cdots=q_d=n^{-2\alpha}$, for some constant $0<\alpha<1$, in this case we have $q\asymp_dq_1^d=n^{-2d\alpha}$. Applying Theorem~\ref{thm19} with $p(1)=1/2$, $n'_1=\cdots=n'_d=n^{2\alpha}$, we obtain the bound $|K|\ll_d n^{d\alpha}
$ under the conditions of Theorem~\ref{tnv2} with $n'=2\,d\,n'_1$. Theorem~\ref{tnv2} itself gives in this case the bound $|K|\ll_d n^{(2d-1)\alpha}$ only (which is worse for $d>1$).

Clearly, Theorem~\ref{thm19} can provide stronger bounds than Theorem~\ref{tnv2} even if relation~\eqref{n11sp} is replaced by some weaker conditions, for example, if
 \begin{equation}
q\ll_d \Big(\prod_{j=1}^d q_j\Big)^{\beta}, \quad \hbox{for some }\beta\ge1/2.
 \label{12sps}
\end{equation}
Conditions~\eqref{11sps} or~\eqref{12sps} are usually satisfied for non-degenerate distributions~$F_a$, for example, if $F_a$ is close to a non-degenerate Gaussian distribution.

Theorem~\ref{tnv2} will be stronger than Theorem~\ref{thm19} if
\begin{equation}
  q\asymp_d q_1=\cdots=q_d.
 \label{22s}
\end{equation}
This is possible, for instance, if the distribution~$F_a$ is close to almost degenerate one which is concentrated in a neighborhood of a one-dimensional subspace. Degenerate distributions can turn, however, into non-degenerate ones after applying linear operators.

Note that we could derive a multivariate analogue of
Theorem~\ref{tnv2} from its one-dimensional version arguing
precisely as in the proof of our Theorem \ref{nthm8}. Then we get
inequality \eqref{n11sp} instead of~\eqref{12sp}.

\bigskip

Similarly as in \cite{EGZ2}, now we state analogues of Theorems~\ref{thm17}--\ref{thm19}
for GAPs of logarithmic rank and with special dimensions, all equal to 1.
Theorems~\ref{tnv3}--\ref{thm55} below extend Theorems~5--7 in~\cite{EGZ2}.

\begin{remark}\label{rem9}\rm In Theorems~\ref{tnv3}--\ref{thm55}, we use the convention 0/0=1.
\end{remark}

\begin{theorem}\label{tnv3}
     Let the conditions of Theorem $\ref{thm17}$ be satisfied,
     except for conditions~\eqref{2b1} and\/~\eqref{2b13}.
Then there exist an absolute positive constant $c_8$ and a GAP $P$ of rank~$r\in\mathbf N$, of volume~$3^r$,
with generators $g_k\in\mathbf{R}$, $k=1,\ldots,r$, and  such that its image $K\subset\mathbf{R}$ has the form
\begin{equation}
{K}=\Big\{\sum_{k=1}^r m_k\, g_k:m_k\in \{-1,0,1\}, \hbox{
for }k=1,\ldots,r\Big\}.\label{1s17}
\end{equation}
Moreover, in the case $\tau>0$ we have
\begin{equation}
r\le c_8\,\big(\left|\log
q\right|+\log(\varkappa/\delta)+1\bigr), \label{n1ss65}
\end{equation}
and at least $n-n'$ elements of\/ $a$ are  $\delta$-close to
 $K$, where $n'\in\mathbf N$ and
\begin{equation}
n'\le c_8\,\big(p(\tau/\varkappa)\big)^{-1}
\bigl(\left|\log q\right|+\log(\varkappa/\delta)+1\bigr)^3.
\label{n1ss68d}
\end{equation}

 In the case $\tau=0$ we have
\begin{equation}
r\le c_8\,(\left|\log
q\right|+1\bigr), \label{n1t6}
\end{equation}
and at least $n-n'$ elements of\/ $a$  belong to
 $K$, where
\begin{equation}
n'\le c_8\,\big(p(0)\big)^{-1}
\bigl(\left|\log q\right|+1\bigr)^3.
\label{n1t8}
\end{equation}
\end{theorem}

\begin{theorem}\label{tnv4}
     Let the conditions of Theorem $\ref{thm27}$ be satisfied, except
     for conditions~\eqref{2s45},~\eqref{2s46}.
Below\/ $c_8$ denotes the
absolute positive constant from Theorem~$\ref{tnv3}$.
     Then, for each $j=1,\ldots,d$,  there exists a GAP $P_j\subset\mathbf{R}$
     of rank~$r_j\in\mathbf N$, of volume~$3^{r_j}$,
with generators $g_k^{(j)}\in\mathbf{R}$, $k=1,\ldots,r_j$, and
such that its image $K_j\subset\mathbf{R}$ has the form
\begin{equation}
{K_j}=\Big\{\sum_{k=1}^{r_j} m_k\, g_k^{(j)}:m_k\in \{-1,0,1\}, \hbox{
for }k=1,\ldots,r_j\Big\}.\label{1s174}
\end{equation}
Moreover, in the case $\tau_j>0$ we have
\begin{equation}
r_j\le c_8\,\big(\left|\log
q_j\right|+\log(\varkappa_j/\delta_j)+1\bigr), \label{n1ss654}
\end{equation}
and at least $n-n_j'$ elements of\/ $a^{(j)}$ are  $\delta_j$-close to
 $K_j$, where $n'_j\in\mathbf N$ satisfy
\begin{equation}
n'_j\le c_8\,\bigl(p(\tau_j/\varkappa_j)\bigr)^{-1}
\bigl(\left|\log q_j\right|+\log(\varkappa_j/\delta_j)+1\bigr)^3.
\label{n1ss68d4}
\end{equation}

 In the case $\tau_j=0$ we have
\begin{equation}
r_j\le c_8\,(\left|\log
q_j\right|+1\bigr), \label{n1t64}
\end{equation}
and at least $n-n'_j$ elements of\/ $a^{(j)}$ belong to
 $K_j$, where
\begin{equation}
n'_j\le c_8\,\bigl(p(0)\bigr)^{-1}
\bigl(\left|\log q_j\right|+1\bigr)^3.
\label{n1t84}
\end{equation}

Define $K=\times
_{j=1}^{d}K_j$. Then the set $K$ is the image of the $d$-dimensional GAP $P$ with rank
\begin{equation}
R=\sum_{j=1}^{d}r_j\le c_8\sum_{j=1}^{d}\bigl(\left|\log
q_j\right|+\log(\varkappa_j/\delta_j)+1\bigr), \label{n1ss658}
\end{equation}
and such that at least  $n-\sum_{j=1}^dn'_j$ elements of\/ $a$ belong to
the set $\times
_{j=1}^{d}[K_j]_{\delta_j}$. Here
\begin{equation}
\sum_{j=1}^dn'_j\le c_8\sum_{j=1}^{d}\bigl(p(\tau_j/\varkappa_j)\bigr)^{-1}
\bigl(\left|\log q_j\right|+\log(\varkappa_j/\delta_j)+1\bigr)^3.
\label{n1ss68d8}
\end{equation}

Furthermore, the set $K$  can be
represented as \begin{equation}
{K}=\Big\{\sum_{k=1}^{R} m_s\, g_s:m_s\in \{-1,0,1\}, \hbox{
for }s=1,\ldots,R\Big\}.\label{1s58}
\end{equation} Moreover, every vector $g_s\in
{\mathbf R}^d$, $s=1,\ldots,R$, has one non-zero coordinate only. Denote
$$
 s_0=0\quad\hbox{and}\quad s_j=\sum_{m=1}^{j}r_m, \quad j=1,\ldots,d.
$$
For $s_{j-1}<s\le s_j$, the vectors $g_s$ are non-zero in the
$j$-th coordinates only and these coordinates are equal to the elements of the
sequence $g_1^{(j)},\ldots, g_{r_j}^{(j)}$ from \eqref{1s174}.
\end{theorem}

\begin{theorem}\label{thm5}  Let\/  $A>0$ and\/ $B\ge0$
be constants. Let\/ $b_n,\varkappa_n, \delta_n,
\tau_n>0$, $n=1,2,\ldots$, be depending on $n$ non-random
parameters satisfying the relations  \;$
b_n^{-B}\leq\delta_n/\varkappa_n\leq1$,
$\delta_n\le\max\{\varkappa_n, \tau_n\}$, for all $n\in \mathbf
N$, and\/ $b_n\to\infty$ as $n\to\infty$.
 Let
 $q_j=Q(F_a^{(j)}, \tau_n)\ge b_n^{-A}$, for  $j=1,\ldots,d$.
 Then, for each $j=1,\ldots,d$,  there exists a GAP $P_j\subset\mathbf{R}$
 of rank~$r_j$, of volume~$3^{r_j}$,
with generators $g_k^{(j)}\in\mathbf{R}$, $k=1,\ldots,r_j$, and
such that its image $K_j\subset\mathbf{R}$ has the form
\begin{equation}
{K_j}=\Big\{\sum_{k=1}^{r_j} m_k\, g_k^{(j)}:m_k\in \{-1,0,1\}, \hbox{
for }k=1,\ldots,r_j\Big\}.\label{1s19}
\end{equation} Moreover, the set $K=\times
_{j=1}^{d}K_j$ is the image of the $d$-dimensional GAP $P$ with rank
\begin{equation}
R=\sum_{j=1}^{d}r_j\ll d\,\big((A+B)\,\log b_n+1\big), \label{21st65}
\end{equation}
and such that at least  $n-n'$ elements of\/ $a$ belong to
the set $\times
_{j=1}^{d}[K_j]_{\delta_n}$. Here $n'\in\mathbf N$ and
\begin{equation}
n'\ll d \,\bigl(p(\tau_n/\varkappa_n)\bigr)^{-1}\big((A+B)\,\log
b_n+1\big)^3. \label{21st68d}
\end{equation}
Furthermore, the description of the set $K$ at the end of the formulation of
Theorem~$\ref{tnv4}$ remains true.\end{theorem}

\begin{theorem}\label{thm55}
The statements of Theorems\/ $\ref{tnv4}$ and\/ $\ref{thm5}$ hold when replacing
$p(\tau_j/\varkappa_j)$ or $p(\tau_n/\varkappa_n)$ by $p (0)$
 in the particular case, where the
parameters $\tau_j$, $j=1,\dots,d$, or $\tau_n$, $n\in\mathbf N$, involved in the formulations
of these theorems, are all zero.
\end{theorem}\medskip

In Theorems 5--7 in~\cite{EGZ2}, we obtained particular cases
of our Theorems~\ref{tnv3}--\ref{thm55}, where
$b_n=n$ and $\tau=\varkappa$, $\tau_j=\varkappa_j$, $j=1,\dots,d$, or $\tau_n=\varkappa_n$, $n\in\mathbf N$.

In Theorems~\ref{tnv3}--\ref{thm55} the approximating GAP may be non-proper.
We could try to get the results with proper GAPs as in Theorems~\ref{thm17}--\ref{thm19}, but then we will lose
the nice representations for the images of GAPs, see \eqref{1s17}, \eqref{1s174}, \eqref{1s19}.
Moreover, the ranks of GAPs will be too large.
\bigskip

\section{Proof of Theorem \ref{thm17}}

 \noindent {\it Proof of Theorem\/ $\ref{thm17}$.}
 Applying  Theorem \ref{thm7} with
 $\delta>0$,
   we derive that, for
$r\in{\mathbf N}_0$, $m\in{\mathbf N}$,
\begin{equation} Q(F_a, \tau)\le c_4^{r+1}\,\frac{\varkappa}{\delta}\,
\biggl(\frac{1}{m\sqrt{p(\tau/\varkappa)\,\beta_{r,m}(M^*,
\delta)}}
+\frac{(r+1)^{5r/2}}{(p(\tau/\varkappa)\,\beta_{r,m}(M^*,
\delta))^{(r+1)/2}}\biggr), \quad\hbox{if }\tau>0, \label{n1sy44}
\end{equation}
and
\begin{equation} Q(F_a, \tau)\le c_4^{r+1}\,
\biggl(\frac{1}{m\sqrt{p(0)\,\beta_{r,m}(M^*, \delta)}}
+\frac{(r+1)^{5r/2}}{(p(0)\,\beta_{r,m}(M^*,
\delta))^{(r+1)/2}}\biggr), \quad\hbox{if }\tau=0, \label{k44}
\end{equation}
with $c_4=2\,c_3$, where $c_3$ is the constant from Theorem~\ref{thm7}. We
assert that the $c_4$ in \eqref{n1sy44}--\eqref{k44} may be taken as the $c_4$ in Theorem
\ref{thm17}.

Let $r\in{\mathbf N}_0$ be fixed and $\tau>0$. Choose now a
positive integer $m=\lfloor y\rfloor+1$, where
 \begin{equation}
y= \frac{2\,c_4^{r+1}\,\varkappa}{
q\,\delta\,\sqrt{p(\tau/\varkappa)\, n'}}\le m. \label{n2s44}
\end{equation}
Assume that $\beta_{r,m}(M^*, \delta)>n'$. Recall that $n'\ge
\big(\,2\,c_4^{r+1}\,(r+1)^{5r/2}\,\varkappa\big/q\,\delta\,\big)^{2/(r+1)}\big/p(\tau/\varkappa)$.
Then, using \eqref{n1sy44} and our assumptions, we have
\begin{equation}
q <q/2+q/2=q. \label{n2s55}
\end{equation}
This leads to a contradiction with the assumption
$\beta_{r,m}(M^*, \delta)>n'$.
 Hence we conclude
$\beta_{r,m}(M, \delta)\le\beta_{r,m}(M^*, \delta)\le n'$.

This means that at least $n-n'$ elements of $a$ are $\delta$-close
to a CGAP $K\in\mathcal{K}_{r,m}$ admitting
representation~\eqref{1s1}, where $h$ is a $r$-dimensional vector,
$V$ is a symmetric convex subset of~${\mathbf R}^r$ containing not
more than $m$ points with integer coordinates. Now equality
\eqref{n2s44} implies the inequality
\begin{equation}
 m\le\frac{2\,c_4^{r+1}\,\varkappa}{ q\,\delta\,\sqrt{p(\tau/\varkappa)\,
n'}}+1, \quad\hbox{if }\tau>0. \label{2s67}
\end{equation}

Without loss of generality we can assume that the absolute values
of $a_k$, $k=1, \ldots,n$, are non-increasing:
$\left|a_1\right|\ge\cdots\ge\left|a_n\right|$. Then, it is easy
to see that
$\left|a_n\right|\le\cdots\le\left|a_{n'+1}\right|\le\left\|a\right\|/\sqrt{n'}$.

If $\delta>\left\|a\right\|/\sqrt{n'}$, we can take as $K^*$
the GAP having zero as the unique element. Then $a_k\in
[K^*]_\delta$, $k=n'+1,\ldots,n$.

Let $\delta\le\left\|a\right\|/\sqrt{n'}$ and
\begin{equation}
V^*=V\cap\big\{{x}\in {\mathbf R}^r:\left|\langle{x},
h\rangle\right|\le2\,\left\|a\right\|/\sqrt{n'}\big\} \label{2t7}
\end{equation}
In this case we take
\begin{equation}
K^*=\big\{\langle{\nu}, h\rangle:{\nu}\in {\mathbf Z}^r\cap
V^*\big\}\subset{\mathbf R} .\label{12s1}
\end{equation}
It is clear that $V^*\subset V$ is a symmetric convex subset
of~${\mathbf R}^r$ and $K^*\in\mathcal{K}_{r,m}$. Moreover,
\begin{equation}
K^*=K\cap\big[-2\left\|a\right\|/\sqrt{n'},\;2\left\|a\right\|/\sqrt{n'}\,\big]
.\label{123s1}
\end{equation}
If $a_k\in [K]_\delta$ and
$\left|a_k\right|\le\left\|a\right\|/\sqrt{n'}$, then $a_k\in
[K^*]_\delta$. Thus, only $n'$ elements of $a$, namely,
$a_1,\ldots,a_{n'}$ may be contained in $[K]_\delta$ and not
contained in $[K^*]_\delta$. Therefore, at least $n-2\,n'$
elements of $a$ are $\delta$-close to the CGAP
$K^*\in\mathcal{K}_{r,m}$.

By Lemma \ref{lmj} and Corollary \ref{lmjj}, there exists a
symmetric, infinitely proper GAP $P$ in ${\mathbf Z}^r$  with rank
$\overline l\le r$ such that we have
\begin{equation}
\label{6s22}\hbox{Image}(P) \subset V^*\cap{\mathbf Z}^r\subset
\hbox{Image}( P_0), \quad  P_0= P^{(c_1r)^{3r/2}},
\end{equation} with
absolute constant~$c_1\ge1$ from the statement of Lemma \ref{lmj}
and
\begin{equation} \hbox{size}(
P_0)\le\big(2\,(c_1\,r)^{3r/2}+1\big)^r\left|V^*\cap{\mathbf
Z}^r\right|\le\big(2\,(c_1\,r)^{3r/2}+1\big)^r\left|V\cap{\mathbf
Z}^r\right|\le\big(2\,(c_1\,r)^{3r/2}+1\big)^rm.\label{777s2}
\end{equation}
Moreover, the generators $g_j$ of $ P_0$, for $1 \le j\le
\overline l$, are contained in the symmetric body $\overline l\hskip1pt V^*$.

Let $\phi:{\mathbf R}^r\to{\mathbf R}$ be a linear map defined
by $\phi(y)=\langle y,h \rangle$, where $h\in{\mathbf R}^r$ is
involved in the definition of $K$ and $K^*$. Define now the
symmetric GAP $\overline P$ with $\hbox{Image}(\overline
P)=\overline K=\big\{\phi(y):y\in \hbox{Image}( P_0)\big\}$ and
with generators $\overline g_j=\phi(g_j)$, $j=1,\ldots,\overline
l$. Obviously, $K^*\subset\overline K\subset{\mathbf R}$ and
$\overline P$ is a symmetric GAP of rank $\overline l$ and
size~$\ll (c_5\,r)^{3r^2/2}m$. The generators $\overline g_j$ are contained
in the set $\big\{\phi(y):y\in \overline lV^*\big\}$. Hence, they
satisfy inequality
 $\bigl|\overline g_j\bigr|\le 2\,r
\left\|a\right\|/\sqrt{n'}$ (see \eqref{12s1} and~\eqref{123s1}).

Applying Lemma~\ref{lmj4} to the GAP $\overline P$, we see that,
for any $t \ge1$, there exists a $t$-proper symmetric
one-dimensional GAP $\overline {\overline P}$ with $\hbox{\rm
rank}(\overline {\overline P})\le \hbox{\rm rank}(\overline P)$,
${\rm Image}(\overline P) \subset {\rm Image}(\overline {\overline
P})\subset{\mathbf R}$, and
\begin{equation} \hbox{\rm size}(\overline P)\le\hbox{\rm size}(\overline {\overline P})
\le(2\,t)^{r}r^{6r^2}\hbox{\rm size}(\overline
P)\ll(2\,t)^{r}r^{6r^2}(c_5\,r)^{3r^2/2}m.\label{86s}
\end{equation}
The statement of item 4 follows from~\eqref{86s} if we take $t=1$.

Let now $t=(c_8r)^{3r/2}$, where $c_8$ is a sufficiently large
absolute constant such that $c_8^{3r/2}\ge2\,r^{5/2}c_1^{3r/2}$.
Let $ \overline {\overline P} = (\overline {\overline L},
\overline {\overline g}, k)$, that is,
\begin{multline}\overline{\overline
K}=\hbox {Image}(\overline {\overline P}) = \big\{ m_1\overline
{\overline g}_1 + \cdots + m_{ k}\overline {\overline g}_{
k}:-\overline {\overline L}_j \le m_j \le \overline {\overline
L}_j,
 \ m_j\in{\mathbf Z},\\ \hbox{ for all }j \hbox{ such that }1 \le j\le   k=\hbox{\rm
rank}(\overline {\overline P}) \le r\big\} \end{multline} with
some generators $\overline {\overline g}_j\in {\mathbf R}$.

Let us prove items 5 and 6. If
$\delta>\left\|a\right\|/\sqrt{n'}$, then we can take as $K^{**}$
the GAP having zero as the unique element. Then $a_k\in
[K^{**}]_\delta$, $k=n'+1,\ldots,n$.

Let $\delta\le\left\|a\right\|/\sqrt{n'}$. It is easy to see that
$\hbox {Image}(\overline {\overline P})  \subset\mathbf R$ is the
image of the box $$B = \big\{(m_1/\overline {\overline
L}_1,\ldots, m_k/\overline {\overline L}_k):(m_1,\ldots, m_k)
\in{\mathbf Z}^k: -\overline {\overline L}_j \le m_j \le \overline
{\overline L}_j\ \hbox{ for all }j=1,\ldots,k\big\}$$ under the
linear map
$$
\Phi:(m_1/\overline {\overline L}_1,\ldots, m_k/\overline
{\overline L}_k) \rightarrow  m_1\overline {\overline g}_1 +
\cdots + m_k\overline {\overline g}_k.
$$
Let now \begin{equation}W = \big\{x=(x_1,\ldots, x_k) \in{\mathbf
R}^k: \left|{x}_j \right| \le 1\ \hbox{ for all
}j=1,\ldots,k\big\}\label{ftt}
\end{equation} and
$$W_t = \big\{tx \in{\mathbf R}^k: x\in W\big\}= \big\{x \in{\mathbf R}^k: \left|{x}_j
\right|\le t\ \hbox{ for all }j=1,\ldots,k\big\}.$$ Let
$$
\Lambda=\big\{(m_1/\overline {\overline L}_1,\ldots, m_k/\overline
{\overline L}_k):(m_1,\ldots, m_k) \in{\mathbf Z}^k\big\}.$$
Obviously, $\Lambda$ is a lattice in~${\mathbf R}^k$. Let
$u=(u_1,\ldots,u_k)\in{\mathbf R}^k $, where $u_j=\overline
{\overline L}_j\,\overline {\overline g}_j$, for $j=1,\ldots,k$.
Then
$$\overline{\overline
K}=\hbox {Image}(\overline {\overline P})  = \big\{\langle{\nu},
u\rangle  :\nu\in \Lambda\cap W \big\}.$$

Define
\begin{equation}
W^*=W\cap\big\{{x}\in {\mathbf R}^k:\left|\langle{x},
u\rangle\right|\le2\,\left\|a\right\|/\sqrt{n'}\big\} \label{f2t7}
\end{equation}
Now we define
\begin{equation}
K^{**}=\big\{\langle{\nu}, u\rangle:{\nu}\in \Lambda\cap
W^*\big\}\subset{\mathbf R} .\label{f12s1}
\end{equation}
It is clear that $W^*\subset W$ is a symmetric convex subset
of~${\mathbf R}^k$ and $K^{**}$ is a CGAP. Moreover,
\begin{equation}
K^{**}=\overline{\overline
K}\cap\big[-2\left\|a\right\|/\sqrt{n'},\;2\left\|a\right\|/\sqrt{n'}\,\big]
.\label{f123s1}
\end{equation}
If $a_k\in [\overline{\overline K}]_\delta$ and
$\left|a_k\right|\le\left\|a\right\|/\sqrt{n'}$, then $a_k\in
[K^{**}]_\delta$. Thus, only $n'$ elements of $a$, namely,
$a_1,\ldots,a_{n'}$ may be contained in $[\overline{\overline
K}]_\delta\setminus[K^{**}]_\delta$. Note that
while counting approximated points, we have already taken into
account that these elements may be not approximated. Therefore, at
least $n-2\,n'$ elements of $a$ are $\delta$-close to the CGAP
$K^{**}$.

By Lemma \ref{lmj} and Corollary \ref{lmjj}, there exists a
symmetric, infinitely proper GAP $R= (N, w, \widetilde l)$ in
$\Lambda$  with rank $\widetilde l\le k\le r$ such that we have
\begin{equation}
\label{f6s22}\hbox{Image}(R) \subset\Lambda\cap W^*\subset
\hbox{Image}( R_0), \quad  R_0= R^{(c_1r)^{3r/2}},
\end{equation}
and
\begin{eqnarray} \hbox{size}(
R_0)&\le&\big(2\,(c_1\,r)^{3r/2}+1\big)^r\left|\Lambda\cap
W^*\right|\le \big(2\,(c_1\,r)^{3r/2}+1\big)^r\left|\Lambda\cap
W\right|\nonumber\\ &\le&
\big(2\,(c_1\,r)^{3r/2}+1\big)^r\hbox{\rm size}(\overline
{\overline P})\le\big(2\,(c_1\,r)^{3r/2}+1\big)^r
(2\,t)^{r}r^{6r^2}(c_5\,r)^{3r^2/2}m.\label{f777s2}
\end{eqnarray}
Moreover, the common generators $w_j$ of $R$ and $ R_0$, for $1
\le j\le \widetilde l$, are contained in the symmetric body $\widetilde l
\,W^*\subset\widetilde l \,W\subset r \,W$ and $w_j\in\Lambda$. In
particular, together with \eqref{ftt} this implies that
$\bigl\|w_j\bigr\|\le r^{3/2} $.

The linear map $\Phi:{\mathbf R}^k\to{\mathbf R}$ can be written
as $\Phi(y)=\langle y,u \rangle$, $y\in{\mathbf R}^k$, where
$u\in{\mathbf R}^k$ is defined above. Define now the symmetric GAP
$\widetilde P$ with $\hbox{Image}(\widetilde P)=\widetilde
K=\big\{\Phi(y):y\in \hbox{Image}( R_0)\big\}$ and with generators
$\widetilde g_j=\Phi(w_j)$, $j=1,\ldots,\widetilde l$. Obviously,
$K^{**}\subset\widetilde K\subset{\mathbf R}$ and $\widetilde P$
is a symmetric GAP of rank $\widetilde l$. The generators
$\widetilde g_j$ are contained in the set $\big\{\Phi(y):y\in \widetilde
l\,W^*\big\}$. Hence, they satisfy the inequality
 $\bigl|\widetilde g_j\bigr|\le 2\,r
\left\|a\right\|/\sqrt{n'}$ (see \eqref{f12s1} and
\eqref{f123s1}).

Let $ \nu\in\hbox{Image}(R_0)$ and $\widetilde t=(c_1r)^{3r/2}$.
Then $\nu$ has a unique representation in the form $$ \nu=m_1\,w_1
+ \cdots + m_{\widetilde l}\,w_{ \widetilde l},$$ where
$-\widetilde t\,N_j \le m_j \le \widetilde t\,N_j$,
 $m_j\in{\mathbf Z}$, for all $j$ such that $1 \le j\le   \widetilde l=\hbox{\rm
rank}(R_0) \le r$.  If $N_j\ge1$, then $\lfloor N_j\rfloor\ge1$
and $\lfloor N_j\rfloor\, w_j\in \hbox{Image}(R) \subset W$.
Therefore, $N_j \,\bigl\|w_j\bigr\|\le 2\,\lfloor N_j\rfloor\,
\bigl\|w_j\bigr\|\le2\sqrt r $. Thus, for all $j=1,\ldots,
\widetilde l$, we have $N_j\, \bigl\|w_j\bigr\|\le 2\,r^{3/2} $.
Hence, $ \left\|\nu\right\|\le2\,\widetilde t\,r^{5/2}\le t $ and
$\nu\in \Lambda\cap W_t$. Since the GAP $\overline {\overline P}$
is $t$-proper, all points of the form $\langle{\nu}, u\rangle $, $
\nu\in \Lambda\cap W_t$, are distinct. The same {can} be said
about all points of the form $\langle{\nu}, u\rangle $, $
\nu\in\hbox{Image}(R_0)\subset \Lambda\cap W_t$. This implies that
the GAP ${\widetilde P}$ is proper. The size of the GAP
$\widetilde P$ coincides with that of $R_0$. It is estimated
{by the right-hand size of \eqref{f777s2} which is $\ll
(c_7\,r)^{21r^2/2}m$ with an absolute constant~$c_7$.} {This
completes the proof of Theorem~\ref{thm17} for $\tau>0$.}

 The case $\tau=0$ can be considered
similarly while using  \eqref{k44} instead
of~\eqref{n1sy44}. Theorem~\ref{thm17} is proved. $\square$
 \medskip

\section{Proof of Theorems $\ref{thm16}$ and $\ref{thm19}$}

\noindent {\it Proof of Theorem\/ $\ref{thm16}$.} First we will
prove Theorem~\ref{thm16} for $d=1$. Denote $q=Q(F_a, 0)$.
  Let $r=r(A, \theta)\in\mathbf{N}_0$ be the minimal
non-negative  integer such that $A<\theta\,(r+1)/2 $. Thus,
$r\le2\,A/\theta$ and $ b_n^{A}<b_n^{\theta(r+1)/2}$ for
all~$b_n>1$. Recall that $b_n\to\infty$ as $n\to\infty$. Assume
without loss of generality that $n$ is so large that
 \begin{eqnarray}\big(2\,c_4^{r+1}\,(r+1)^{5r/2}/q\big)^{2/(r+1)}/p(0)&\le&
\big(2\,c_4^{r+1}\,(r+1)^{5r/2}\varepsilon_1^{-1}\,b_n^{A}\big)^{2/(r+1)}/p(0)\nonumber\\
&\le& \varepsilon_2\,b_n^\theta\le n'=n'_1 \le
n.\label{3uu}\end{eqnarray} It remains to apply
Theorem~\ref{thm17} with $\tau=\delta=0$.

If \eqref{3uu} is not satisfied, then $n'\le n=O(1)$ and  we can
take as $K$ the set
\begin{equation}\label{3uut}K(n')=\Big\{\sum_{k=\,n'+1}^{n} s_k a_k:s_k\in \{-1,0,1\}, \hbox{ for }k=1,\ldots,r\Big\}.
 \end{equation}
Without loss of generality we can assume again that the absolute
values of $\left|a_k\right|$, $k=1, \ldots,n$, are non-increasing:
$\left|a_1\right|\ge\cdots\ge\left|a_n\right|$. Clearly, $K(n')$
is the image of a GAP $P(n')$ of rank $n-n'$ and of
size~$3^{n-n'}$.
 The generators of $P(n')$ are $a_{n'+1},\ldots,a_{n}$ satisfying
 $\left|a_n\right|\le\cdots\le\left|a_{n'+1}\right|\le\left\|a\right\|/\sqrt{n'}$.
At least $n-n'$ elements of $a$ (namely, $a_{n'+1},\ldots,a_{n}$)
belong to $K(n')$. Of course, the gap $P(n')$ may be non-proper.
In order to find a proper gap $P$, one should proceed like as in
the proof of Theorem~\ref{thm17}, using Lemma \ref{lmj} and
Corollary \ref{lmjj}. Thus Theorem~\ref{thm16} is proved
 for $d=1$. \medskip

Let us now assume that $d>1$. We apply Theorem \ref{thm16}
with $d=1$ to the distributions of the coordinates of the
vector~$S_a$,
 taking the vector $a^{(j)}=(a_{1j},\ldots,a_{nj})$ as vector $a$, $j=1,\ldots, d$.
 Then, for each $a^{(j)}$, there exists a proper symmetric GAP $P_j$ with image $K_j \in \mathcal{K}_{r_j,m_j}$,
 which satisfies the assertion of Theorem~\ref{thm16}, that is:

1. At least $n-2\,n'_j$ elements of $a^{(j)}$ are contained in
 $K_j$;

2. $P_j$ has small rank $l_j\le r_j=O(1)$, and
 \begin{equation}
 m_j\le\max\Big\{O\Big(q_j^{-1}\,(n'_j)^{-1/2}\Big),
1\Big\}. \label{2s44788}
\end{equation}

3. The generators ${g}_1^{(j)},\ldots, {g}_{l_j}^{(j)}$ of~$P_j$ satisfy the inequality
 $\bigl|g_p^{(j)}\bigr|\le
2\,r_j\, \left\|a^{(j)}\right\|/\sqrt{n'_j}$, for \;$p=1,\ldots,l_j$.

Thus, the multi-vector $a$ is well approximated by
 the GAP $P$ with image $K=\times _{j=1}^{d} K_j$ of rank $l_1+\cdots +l_d=L\le R=r_1+\cdots +r_d$.
 At least $n-2\sum_{j=1}^dn'_j$ elements of\/ $a$ are contained in~$K$.

 It is easy to see that $K\in\mathcal{K}_{r,m}^{(d)} $,
 $r=(r_1, \ldots,r_d)\in{\mathbf N}_0^d$, $m=(m_1, \ldots,m_d)\in{\mathbf N}^d$,
\begin{equation}\label{32dd4}\left|\times _{j=1}^{d}K_j\right|\le \prod_{j=1}^{d} m_j\le\prod_{j=1}^{d}
 \max\Big\{O\Big(q_j^{-1}\,(n'_j)^{-1/2}\Big),
1\Big\}.\end{equation}
\hfill $\square$
\medskip

 \noindent {\it Proof of Theorem\/ $\ref{thm19}$.} First we will
 prove Theorem~\ref{thm19} for $d=1$. Denote $q=Q(F_a, \tau_n)$.
   Let $r=r(A, B,\theta,D)$ be the minimal
positive  integer such that $A+B<(\theta-D)\,(r+1)/2 $. Thus,
$r\le2\,(A+B)/(\theta-D)$ and $b_n^{A+B}< b_n^{(\theta-D)(r+1)/2}$
for all~$b_n>1$. Recall that $b_n\to\infty$ as $n\to\infty$.
Assume without loss of generality that $n$ is  large enough such that
\begin{eqnarray}\big(2\,c_4^{r+1}\,(r+1)^{5r/2}\,\varkappa_n/q\,\delta_n\,\big)^{2/(r+1)}/p(\tau_n/\varkappa_n)
&\le&
\big(2\,c_4^{r+1}\,(r+1)^{5r/2}\varepsilon_1^{-1}\,
\varepsilon_4^{-1}\,b_n^{A+B}\big)^{2/(r+1)}/\varepsilon_3\,b_n^{-D}\nonumber\\
&\le& \varepsilon_2\,b_n^\theta\le n' \le
n.\label{03dd}\end{eqnarray} It remains to apply
Theorem~\ref{thm17}. If \eqref{03dd} is not satisfied, then
$n=O(1)$ and  we can again take as $K$ the set $K(n')$ defined in
\eqref{3uut}. In order to find a proper gap $P$, one should
proceed  as in the proof of Theorem~\ref{thm17}.
Thus Theorem~\ref{thm19} is proved
 for $d=1$. \medskip

Let us now assume that $d>1$. We apply Theorem \ref{thm19} with $d=1$ to the
distributions of the coordinates of the vector~$S_a$,
 taking the vector $a^{(j)}=(a_{1j},\ldots,a_{nj})$ as  vector $a$, $j=1,\ldots, d$.
 Then, for each $a^{(j)}$, there exists a  proper symmetric GAP $P_j$ with image $K_j \in \mathcal{K}_{r_j,m_j}$,
 which satisfies the assertion of Theorem~\ref{thm19}, that is:

1. At least $n-2\,n'_j$ elements of $a^{(j)}$ are
$\delta_n$-close to
 $K_j$;

2. $P_j$ has small rank $l_j\le r_j=O(1)$, and
 \begin{equation}
 m_j\le\max\Big\{O\Big(q_j^{-1}\,\rho_n^{-1}\,\big(n'_j\,p(\tau_n/\varkappa_n)\big)^{-1/2}\Big),
1\Big\}. \label{02s4478}
\end{equation}

3. The generators ${g}_1^{(j)},\ldots, {g}_{l_j}^{(j)}$ of~$P_j$,  satisfy the inequality
 $\bigl|g_p^{(j)}\bigr|\le
2\,r_j\, \left\|a^{(j)}\right\|/\sqrt{n'_j}$, for \;$p=1,\ldots,l_j$.

Thus, the multi-vector $a$ is well approximated by
 the GAP $K=\times _{j=1}^{d} K_j$.
 It is easy to see that $K\in\mathcal{K}_{r,m}^{(d)} $,
 $r=(r_1, \ldots,r_d)\in{\mathbf N}_0^d$, $m=(m_1, \ldots,m_d)\in{\mathbf N}^d$,
\begin{equation}\label{032dd}\left|\times _{j=1}^{d}{K_j}\right|\le \prod_{j=1}^{d} m_j\le\prod_{j=1}^{d}
 \max\Big\{O\Big(q_j^{-1}\,\rho_n^{-1}\,\big(n'_j\,p(\tau_n/\varkappa_n)\big)^{-1/2}\Big),
1\Big\}.\end{equation}

Since at most $2\,n'_j$ elements of $a^{(j)}$ are far from the
GAPs $K_j$,
 there are at least $n-2\sum_{j=1}^dn'_j$ elements of $a$ that are  $\delta_n$-close
to the GAP $K$. In view of relation \eqref{032dd} and taking into
account that $K=\times _{j=1}^{d} K_j$, we obtain
relation~\eqref{n111sp}. $\square$\medskip

\begin{remark}\rm Notice that, in Theorems  \ref{thm16} and \ref{thm19}, the
ranks of $P_j$ are actually the same for all $j=1,\ldots, d$.
Moreover, in Theorems  \ref{thm16} and \ref{thm19}, we get explicit bounds for~$r_j$, for
sufficiently large~$n$, namely:
$r_j\le2\,A/\theta$ and $r_j\le2\,(A+B)/(\theta-D)$ respectively.
\end{remark}
\bigskip

\section{Proofs of Theorems \ref{tnv3}--\ref{thm55}}

\noindent {\it Proof of Theorem\/ $\ref{tnv3}$.}
By Corollary~$\ref{lm429}$, we have
\begin{equation}\label{1155t}
q\ll \frac \varkappa\delta \,Q,\quad\hbox{where}\ Q=Q(H^{p(\tau/\varkappa)}, \delta).
\end{equation}
Note that $H^{p(\tau/\varkappa)}$ is the
symmetric infinitely divisible distribution  with L\'evy
 spectral measure $\frac{\,p(\tau/\varkappa)\,}4\;M^*$, where
 $M^*=\sum_{k=1}^{n}\big(E_{a_k}+E_{-a_k}\big)$.

Applying
Theorem 3.3 of Chapter~II
\cite{Arak and Zaitsev} (see Theorem 4 in \cite{EGZ2}), we obtain that
 there exist $r\in\mathbf N_0$ and $g_j\in
{\mathbf R}$, $j=1,\ldots,r$, such that
\begin{equation}
r\ll\left|\log Q\right|+1, \label{pss65}
\end{equation}
and
\begin{equation}
p(\tau/\varkappa) \,M^*\{{\mathbf R}\setminus[K]_{\delta}\}\ll
\bigl(\left|\log Q\right|+1\bigr)^3, \label{pss68d}
\end{equation}
where $K$  has the form \eqref{1s17}. Recall that $\delta\le\varkappa$.
By \eqref{1155t},
\begin{equation}
\left|\log Q\right|\ll \left|\log q\right|+\log\bigl(\varkappa/\delta\bigr). \label{ps6}
\end{equation}
Inequalities \eqref{pss65}--\eqref{ps6} together imply the statement of Theorem\/ $\ref{tnv3}$ in the case $\tau>0$.
The case $\tau=0$ is a little bit easier.  $\square$
\medskip

Theorems~\ref{tnv4}--\ref{thm55} are direct consequences of Theorem~\ref{tnv3}.
\bigskip

We are grateful to anonymous reviewers for useful remarks.

\end{document}